\documentclass[12pt]{amsart}

\usepackage{amsmath}
\usepackage{amscd,amssymb}
\usepackage{epsfig}

\usepackage{graphicx}
\usepackage{color}
\usepackage{hyperref}
\usepackage{booktabs}
\usepackage{amsfonts}
\usepackage{mathrsfs}

\def\qed{\ifmmode\square\else\nolinebreak\hfill
$\Box$\fi\par\vskip12pt}

\newtheorem{thm}{Theorem}[section]
\newtheorem{lemma}[thm]{Lemma}
\newtheorem{corollary}[thm]{Corollary}

\numberwithin{equation}{section}
\numberwithin{thm}{section}

\theoremstyle{definition}
\newtheorem{definition}[thm]{Definition}
\newtheorem{remark}[thm]{Remark}

\newcommand{\bQ}{\mathbb Q}

\newcommand{\bZ}{\mathbb Z}

\newcommand{\cD}{\mathcal D}

\newcommand{\cR}{\mathcal R}

\definecolor{Purple}{rgb}{0.5,0,0.5}

\begin{document}
\pagestyle{plain}
\begin{titlepage}

\title{Two Kinds of Constructions of Directed Strongly Regular Graphs from Partial Sum Families and Semi-direct Products of Groups}
\begin{center}
\author{Jingkun Zhou, Zhiwen He$^*$ and Zhao Chai
}\end{center}
\address{School of Mathematical Sciences, Zhejiang University, Hangzhou 310027, China}
\email{jingkunz@zju.edu.cn}

\address{School of Mathematical Sciences, Zhejiang University, Hangzhou 310027,  China}
\email{zhiwen$\_$he@zju.edu.cn}

\address{School of Mathematical Sciences, Peking University, Beijing 100871, China}

\email{chaizhao@scichina.org}

\begin{abstract} In this paper, we construct directed strongly regular graphs with new parameters by using partial sum families with local rings. 16 families of new directed strongly regular graphs are obtained and the uniform partial sum families are given. Based on the cyclotomic numbers of finite fields, we present two infinite families of directed strongly regular Cayley graphs from semi-direct products of groups.
\end{abstract}

\keywords{directed strongly regular graph, partial sum family, semidirect product, cyclotomic number\\
{\bf  Mathematics Subject Classification (2010) 05C20 05E15 05E18 05E40}\\
{\bf  Funding information: National Natural Science Foundation of China under Grant No. 11771392.}\\
$^*$Correspondence author}

\maketitle


\section{Introduction}
Directed strongly regular graphs(or DSRG, for short) were introduced by Duval \cite{2} in 1988. A DSRG $\mathcal{G}=(V,E)$ with parameters $(v,k,\lambda,\mu,t)$ is a directed graph without loops satisfing that
\begin{itemize}
	\item[(1)] for any vertex $x$, there are exactly $k$ distinct vertices $y$ and $k$ distinct vertices $z$ such that $(x,y),(z,x)\in E$;
	\item[(2)] for any vertex $x$, there are exactly $t$ vertices $y$ such that $(x,y),(y,x)\in E$;
	\item[(3)] for any two distinct vertices $x$ and $y$, the number of vertices $z$ such that $(x,z),(z,y)\in E$ is exactly $\lambda$ if $(x,y)\in E$, or $\mu$ otherwise.
\end{itemize}

For any arc $(x,y)\in E$, if $E$ also contains an arc $(y,x)$, then $\{(x,y),(y,x)\}$ can be viewed as an undirected edge. It is clear that each vertex of $\mathcal{G}$ is on $t$ undirected edges and $2k-2t$ directed edges. Observe that a DSRG is a (undirected) strongly regular graph if $t=k$, or a doubly regular tournament if $t=0$. The two cases, together with the one that $\mathcal{G}$ is complete, are called trivial and we only consider nontrivial (or in other words, genuine) DSRGs in this paper.

Several necessary conditions on the parameters of a DSRG are given in \cite{2}. In addition to some nonexistence proofs, Duval gave the existence results in his paper through block constructions and Kronecker product constructions. The work of Duval has been followed by many researchers and a plenty of new tools and methods have been introduced to construct DSRGs. With the aid of a computer, Fiedler, Klin and Muzychuk \cite{13} determined all DSRGs of order $v\le 20$ which have a vertex-transitive automorphism group. Together with a result by J{\o}rgensen \cite{12}, they gave a complete answer on Duval's question about the existence of directed strongly regular graphs with $v\le20$.  In the light of algebraic characterization of strongly regular graphs, DSRGs were studied by coherent algebras; see \cite{13} and \cite{14}-\cite{16} for details. DSRGs also can be obtained from many other combinatorial structures, including finite incidence structures in \cite{17}-\cite{19}, equitable partitions in \cite{20} and block matrices in \cite{22}.

The Cayley digraph is one of the most important tools to construct SRGs and DSRGs. A necessary and sufficient condition is given in \cite{21} for a subset of a finite group to be the connection set of a directed strongly regular Cayley graph and obtained an infinite family of DSRGs. A necessary and sufficient condition for a subset of a finite group to be the connection set of directed strongly Cayley graph are  given in \cite{21}.
\begin{definition}
Let $m,n$ be two integers with $m\ge1$ and $n\ge2$, and $\mathcal{G}$ be a digraph of $mn$ vertices. If there is an automorphism group $G$ of $\mathcal{G}$ of order $n$ such that the vertex set of $\mathcal{G}$ can be divided into $m$ orbits of length $n$ by the action of $G$, and the actions of $G$ on each orbit are all regular, then $\mathcal{G}$ is said to be an \textbf{$m$-Cayley digraph}. In particular, when $m=1$, we call $\mathcal{G}$ a Cayley digraph; when $m=2$, $\mathcal{G}$ is a semi-Cayley digraph.
\end{definition}
By \cite{12}, a genuine DSRG cannot be a Cayley digraph of an abelian group, and thus two common approaches to construct genuine DSRGs are derived, one of which is to consider Cayley digraphs of nonabelian groups, while the other is to consider $m$-Cayley digraphs for $m\geq 2$, a generalized version of Cayley digraphs. See \cite{3} for a demonstration of the former approach in which Duval and Iourinski constructed a family of directed strongly regular Cayley graphs of certain semidirect product groups. In this paper we also have two new semidirect product constructions of directed strongly regular Cayley graphs using cyclotomic numbers, which will be stated in Section 4.

A succession of papers has revealed that m-Cayley graphs are useful combinatorial objects for studying SRGs: see paper \cite{23}, \cite{24} for semi-Cayley graphs and \cite{25} for tri-Cayley graphs. Mart\'inez and Araluze \cite{6}, based on the earlier work of Duval and Iourinski \cite{3}, provided the concept of partial sum families that generate $m$-Cayley digraphs that are DSRGs.  Also paper \cite{5} gives the formal definition and some new constructions of partial sum families, which confirms the practicability of this object. In \cite{7}, the authors considered the particular case that the group is abelian and $m = 2$, in which the partial sum families are called bi-Abelian partial sum quadruples (or PSQs briefly) and they characterized the parameters of bicirculant PSQs. In the end, they concluded with three open problems, of which the second one is to construct directed strongly regular $m$-Cayley digraphs for each $m\geq 2$. We prove the existence of partial sum families with $m\geq 2$, using the knowledge of group ring and character theory, and give an answer to this problem. Inspired by the method of constructing partial difference sets with Galois rings in \cite{8}, we try to shift our attention to construct partial sum families with local rings.

The complementary graph $\overline{\mathcal{G}}=(V',E')$ of a directed graph $\mathcal{G}=(V,E)$ is a directed graph with vertex set $V'=V$ and arc set $E'=\{(x,y)\in V\times V: (x,y)\notin E\}$. Duval declared in \cite{2} that the complement of a DSRG with parameters $(v,k,\lambda,\mu,t)$ is also a DSRG with parameters $(v',k',\lambda',\mu',t')=(v,v-k-1,v-2k+\mu-2,v-2k+\lambda,v-2k+t-1)$. For any DSRG $\mathcal{G}_1$ with parameters $(v,k,\lambda,\mu,t)$, we call a DSRG $\mathcal{G}_2$ has complementary parameters with respect to $\mathcal{G}_1$ if $\mathcal{G}_2$ is of parameters $(v',k',\lambda',\mu',t')$. In this sense we tend to focus only on DSRGs with parameters such that $k\le v/2$. We compare the parameters which are given by our first construction with the ones listed in Brouwer's homepage which are still unknown, and conclude that there are 16 new 5-tuples of parameters from our construction in the case that the vertex number $v<110$, which are
$$\begin{matrix}
(50, 18, 7, 6, 12), &(75, 28, 11, 10, 16), &(75, 32, 13, 14, 20), &(98, 26, 9, 6, 16),\\
(98, 32, 11, 10, 22), &(98, 39, 16, 15, 27), &(100, 38, 15, 14, 20), &(100, 42, 17, 18, 24)\\
\end{matrix}$$
and their complementary tuples.

The rest of this paper is organized as follows. In Section 2, we give some preliminaries on local rings, $m$-Cayley digraphs, partial sum families, semidirect products and cyclotomic numbers. In Section 3, we obtain, by using local rings, an infinite family of partial sum families, some of which generate DSRGs with new parameters. We also obtain uniform partial sum families under the same frame. In Section 4, we use cyclotomic numbers of finite fields to give two constructions of directed strongly regular Cayley graphs over the semidirect product of the additive group of finite fields. At last, Section 5 concludes the paper.

\end{titlepage}
\section{Preliminaries}

\subsection{Group rings}
Let $G$ be a finite group. The group ring $\bZ[G]$ is defined as the set of the formal sums of elements of $G$ with coefficients in $\bZ$. The operations ``$+$" and ``$\cdot$" of $\bZ[G]$ are given by
\begin{displaymath}
\sum_{g\in G}a_gg+\sum_{g\in G}b_gg=\sum_{g\in G}(a_g+b_g)g
\end{displaymath}
and
\begin{displaymath}
\left(\sum_{g\in G}a_gg\right)\cdot\left(\sum_{h\in G}b_hh\right)=\sum_{g,h\in G}a_gb_h(gh).
\end{displaymath}

Many combinatorial structures were studied within the context of a group ring. It becomes conventional to abuse the notation $S$ as a subset of $G$ and the corresponding element $\sum_{s\in S}s$ in $\bZ[G]$ at the same time. Character theory is another powerful tool which will be used in the following context. The importance of character theory lies that it can help us simplify calculations. A character of a finite abelian group $G$ is a homomorphism from $G$ to the multiplicative group of complex numbers of absolute value 1. The principal character is the character $\chi_0$ such that $\chi_0(x)=1$ for any $x\in G$. We write $\chi(S)=\sum_{s\in S}\chi(s)$ for a subset $S$ of $G$.

We also need the well-known inverse formula on group rings:
{\lemma\label{lem_4}
Suppose $A=\sum_{g\in G}a_g g$ is an element of the group ring $\bZ[G]$ for a finite abelian group $G$, then the coefficients $a_g$'s of $A$ can be computed explicitly by
$$a_g=\frac{1}{|G|}\sum_{\chi\in\hat{G}}\chi(A)\chi(g^{-1}),$$
where $\hat{G}$ denotes the character group of $G$. In particular, if $A,B\in\mathbb{Z}[G]$ satisfy $\chi(A)=\chi(B)$ for all characters $\chi\in\hat{G}$, then $A=B$.}

\subsection{Local rings}

A local ring is a ring which has a unique maximal ideal. In this paper we only consider finite commutative local rings whose maximal ideals are principal. Such local rings are sometimes called finite chain rings. Denote a finite chain ring by $\cR$ and its maximal ideal by $I$, then $I$ is generated by a prime element $\pi$ and consists of all zero divisors of $\cR$ with 0. $\cR\setminus I$ is the set of units. Leung and Ma \cite{4} constructed partial difference sets in $\cR\times\cR$ with a finite chain ring $\cR$ that satisfies some certain properties. A natural idea arises that we can construct partial sum families over the same group $\cR\times \cR$. Here we list some properties and facts appeared in \cite{4} in which will be used in Section 3.
{\lemma\label{lem_1}
(Propsition 2.4, \cite{4}) For any prime integer $p$ and positive integers $s,r,d$ with $r\le s$, there exists a finite chain ring $\cR$ with its maximal ideal $I=(\pi)$ for a prime element $\pi$ such that (i)\ $I^{s-1}\ne0$, $I^s=0$, (ii)\ $p$ is an associate of $\pi^r$ and (iii)\ $\cR/I$ is a finite field with $p^d$ element.
}

People can obtain many finite chain rings which satisfy the conditions stated in Lemma \ref{lem_1}. For example, to our best known, one can choose a finite field or a Galois ring.  A sophisticated example is also readily available in \cite{4}: let $\bQ_p$ be the $p$-adic completion of $\bQ$ and $\bZ_p$ be the set of $p$-adic integers. Let $K$ be a finite extension of $\bQ$ and $R'$ be the integral closure of $\bZ_p$ in $K$. Now $R'$ is also a local ring with maximal ideal $(\pi)$. Then we get the finite chain ring $R:=R'/\pi^sR'$ we wanted.

Let $H$ be the additive group of the ring $\cR$. We consider the additive character group $\hat{H}$ of $H$. Since $I^{s-1}\ne0$, there exists a character $\psi\in\hat{H}$ such that $\psi$ is nonprincipal on $I^{s-1}$. For each $a\in\cR$, the map $\psi_{a}$ defined by $\psi_{a}(x)=\psi(ax)$ for $x\in\cR$ is an additive character of $\cR$. It is clear that $\psi_{a}\ne\psi_{b}$ holds whenever $a\ne b$ since $\psi$ is nonprincipal on $I^{s-1}$. As a result, the set $\{\psi_{a}|a\in\cR\}$ gives all additive characters of $\cR$. We rewrite the facts mentioned above in the following lemma.
{\lemma\label{lem_2}
(\cite{4}) There is an additive character $\psi\in\hat{H}$ which is nonprincipal on $I^{s-1}$. Besides, $\hat{H}=\{\psi_{a}|a\in\cR\}$ where $\psi_a(x)=\psi(ax)$, $x\in\cR$.
}

Now we focus on the group $G=H\times H$ where $H$ is the additive group of a finite chain ring $\cR$ which satisfies the conditions in Lemma 2.1 for positive integers $s,r,d$ and prime integer $p$. There is a simple but useful fact on the description of characters of $G$.
{\lemma\label{lem_char}
(Propsition 3.2, \cite{4}) For any character $\chi\in\hat{G}$, there exist $a,b\in H$ such that $\chi(x,y)=\psi(ax+by)$ for $x,y\in H$.}

In \cite{8} Pohill constructed partial difference sets by forming a ``spread" of the group $G$ when $H$ is a Galois ring. In this paper we follow this idea to construct partial sum families in $G$. We generalize the ``spread" over a Galois ring to the one over a finite chain ring $\cR$:
$$\begin{aligned}
&L_{a}=\{(x,ax)|x\in\cR\}\\
&L_{\infty}=\{(0,x)|x\in\cR\}
\end{aligned}$$
Fix a system of coset representatives of $R$ module $I$ and denote it by $J$.Let $J'=J\cup\{\infty\}$, then for any two distinct elements $a$ and $b$ of $J$, we have $L_a\cap L_b=\{(0,0)\}$. Since $L_a\cap L_b=\{(0,0)\}$ and $|L_a|=|L_b|=|H|$ for any $a\neq b\in J'$, we have $L_aL_b=G$.

\subsection{m-Cayley digraphs and partial sum families}
 Let us briefly describe the notion of a difference digraph and determine when they are DSRGs by using the concept of partial sum families. Assume that $G$ is an abelian group of order $n$ whose operation is additively written. Now define the vertex set $V$ of a digraph $\mathcal{G}$ to be the union of $m$ copies of $G$, that is,
$$V=\bigcup_{0\le i\le m-1}G_{i}\textnormal, $$
where each $G_{i}$ is a copy of $G$. Then there is a left translation action $\rho$ of $G$ on the vertex set $V$ and each $G_{i}$ is an orbit of $\rho$.

Next we define the arc set $E$ of $\mathcal{G}$. Let $\{S_{i,j}\}_{0\le i,j\le m-1}$ be $m^2$ subsets of $G$. For any two distinct vertices $x\in G_i$ and $y\in G_j$, $0\leq i,j\leq m-1$, we set $(x,y)\in E$ if and only if $x-y\in S_{i,j}$. As a result, $G$ is an automorphism group of the digraph $\mathcal{G}=(V,E)$ and acts on each $G_i$ regularly. The digraph $\mathcal{G}$ is called a difference digraph associated with the family $\{S_{i,j}\}_{0\leq i,j\leq m-1}$. Note that if a digraph $\mathcal{G}$ is a difference digraph associated with a family $\{S_{i,j}\}$, then its complement $\overline{\mathcal{G}}$ is a difference digraph associated with the family $\{S'_{i,j}\}$ with
$$S'_{i,i}=G^{\ast}\setminus S_{i,i}$$
for $0\le i\le m-1$ and
$$S'_{i,j}=G\setminus S_{i,j}$$
for $0\le i\ne j\le m-1$, where $G^{\ast}=G\setminus\{e\}$.
In \cite{5}, Araluze, et al. showed that $\mathcal{G}$ is a DSRG if and only if $\{S_{i,j}\}_{0\leq i,j\leq m-1}$ is a partial sum family defined as follows.

\begin{definition}{\label{def_psf}}
 A family $\{S_{i,j}\}_{0\le i,j\le m-1}$ of subsets of a finite abelian group $G$ is called an \textbf{$(m,n,k,\lambda,\mu,t)$-partial sum family}, if members of this family satisfy the following conditions:
\begin{itemize}
	\item[(1)] for any $0\leq i\le m-1$, $e\notin S_{i,i}$, where $e$ is the identity of $G$;
	\item[(2)] for any $0\leq i\le m-1$, $\sum_{j=0}^{m-1}|S_{i,j}|=\sum_{j=0}^{m-1}|S_{j,i}|=k$;
	\item[(3)] for any $0\leq i,j\le m-1$, $\sum_{l=0}^{m-1}S_{l,j}S_{i,l}=\mu G+\beta S_{i,j}+\delta_{i,j}\gamma e$, where $\beta=\lambda-\mu$, $\gamma=t-\mu$, and $\delta_{i,j}$ is the Kronecker delta.
\end{itemize}
\end{definition}

\begin{lemma}{\label{lem_DF}}
(Propsition 1.2, \cite{5}) A family $\{S_{i,j}\}_{0\le i,j\le m-1}$ of subsets of a finite abelian group $G$ is an $(m,n,k,\lambda,\mu,t)$-partial sum family if and only if the difference digraph defined by $\{S_{i,j}\}$ is an $(mn,k,\lambda,\mu,t)$-DSRG.
\end{lemma}
In the case of $m=2$, the regularity of $\mathcal{G}$ gives rise to $|S_{0,0}|=|S_{1,1}|$ and $|S_{0,1}|=|S_{1,0}|$. For $m>2$, however, such desirable properties can not be obtained from the regularity. To require partial sum families for $m>2$ to perform similarly as the ones for $m=2$, the authors gave the notion of uniformity in \cite{7}.
\begin{definition}{\label{def_upsf}}
(\cite{7})
A partial sum family $\{S_{i,j}\}_{0\leq i,j\leq m-1}$ of a finite abelian group $G$ is said to be \textbf{uniform} if it satisfies the conditions:
\noindent(1) $|S_{i,i}|$ are all equal.
\noindent(2) $|S_{i,j}|$, $i\neq j$ are all equal.
\noindent(3) $\{S_{i,i}: 0\leq i\leq m-1\}$ form a partition of $G\setminus\{0\}$.
\end{definition}

\subsection{Semidirect products and Cyclotomic numbers}
In this subsection we introduce cyclotomic numbers over finite fields and semidirect products of groups. Let $\mathbb{F}_{q}$ be a finite field with $q$ elements, where $q$ is a power of an odd prime $p$. Consider the multiplicative group $\mathbb{F}_{q}^{\ast}$. Let $\omega$ be a primitive element of $\mathbb{F}_{q}^{\ast}$, and write $q-1=ef$, where $e,f>1$ are two integers. Let
$$C_{0}=\langle\omega^{e}\rangle=\{\omega^{ej}|0\le j\le f-1\}$$
be the subgroup of $\mathbb{F}_{q}^{\ast}$ generated by $\omega^{e}$, and let
$$C_{i}=\omega^{i}C_{0}=\{\omega^{i}x|x\in C_{0}\},\ \textnormal{for}\ 0\le i\le e-1,$$
be the cyclotomic cosets of $C_{0}$. Then the cyclotomic numbers of order $e$  is defined by
$$(i,j)_{e}=|(C_{i}+1)\cap C_{j}|,$$
for $0\le i,j\le e-1$, and write $(i,j)$ for $(i,j)_{e}$ simply if the value of $e$ is clear from the context.
By \cite{11}, the partition $C_{0},\dots,C_{e-1},\{0\}$ of $\mathbb{F}_{q}$ form a subalgebra, which is called a Schur ring, of the group ring $\mathbb{Z}[H]$, where $H$ is the additive group of $\mathbb{F}_{q}$. In other words, the subspace spanned by $C_{i}$'s and $\{0\}$ is closed under multiplication
$$C_{i}C_{j}=\sum_{k=0}^{e-1}p_{ij}^{k}C_{k}+|C_{i}\cap-C_{j}|[0],$$
with $p_{ij}^{k}\in\mathbb{Z}$ for $0\le i,j,k\le f-1$, where $[0]$ denotes the corresponding element in $\mathbb{Z}[H]$ to 0. In particular,
\begin{equation}\label{Eqn_schurring}
C_{0}C_{i}=\sum_{k=0}^{e-1}(i,k)C_{k}+|C_{0}\cap-C_{i}|[0].
\end{equation}

\begin{definition} Let $H$ and $K$ be two groups and $\theta:K\rightarrow$ Aut$\,H$ be a homomorphism. The \textbf{semidirect product} $H\rtimes K$ of $H$ by $K$ relative to $\theta$ is the group with underlying set $H\times K$ consisting of pairs $(h,k)$ for $h\in H$ and $k\in K$, and operation $$(h,k)(h',k')=(h(\theta(k)(h')),kk').$$
The identity element is $(1_{H},1_{K})$ and the inverse of $(h,k)$ is $(h,k)^{-1}=(\theta(k^{-1})(h^{-1}),k^{-1})$.
\end{definition}
Readers are referred to \cite{26} for more details about semidirect products. We end this section with a lemma given in \cite{21}.

\begin{lemma}
(\cite{21})
The Cayley digraph Cay($G,S$) of a finite group $G$ with respect to a subset $S$ is a DSRG with parameters $(v,k,\lambda,\mu,t)$ if and only if the following conditions hold:

\noindent(1)$|G|=v$; (2)$|S|=k$;

\noindent(3)$S^{2}=te+\lambda S+\mu(G-e-S)=\mu G+\beta S+\gamma e$, where $e$ is the identity element of $G$.
\end{lemma}

\section{Constructions of Directed Strongly Regular Graphs from Partial Sum Families}
Suppose $\cR$ is a local ring of order $p^{sd}$ which satisfies the properties of Lemma 2.1 and $I$ is the unique maximal ideal. Let $H$ be the additive group of $\cR$ and $G=H\times H$.  In this section we give a brief description of our strategy to obtain partial sum families by using the ``spreads'' defined in Section 2.2. Then we get a new infinite family of  DSRGs that derived from these partial sum families by lemma \ref{lem_DF}. When the ring $\mathcal{R}$ has $I=0$, we can also obtain uniform partial sum families. Here we follow the notations $J,J'$ and $L_{a}$'s for $a\in J'$ denoted in Section 2.2. For the identity element $(0,0)$ of $G$, we simply write it by $0_G$.

We first state our main result.

{\thm\label{thm_psfmain}
For any prime integer $p$ and positive integers $s,d,w,z_1$, $1\leq w\leq p^d$, $1\leq z_1\leq p^d-w+1$, there exists a DSRG with parameters $(mn,k,\lambda,\mu,t)=(z_{1}+1,p^{2sd},(p^{sd}-1)w+z_{1}z_{2}p^{sd},p^{sd}+w^2-3w+z_{1}z_{2}^2,w^2-w+z_{1}z_{2}^2,p^{sd}w-w+z_{1}z_{2}^2)$, where integer $z_2=1$ if $w=1$ and $z_2=w-1$ or $w$ if  $2\leq w\leq p^d$.}

Label the elements of $J'$ by $a_0=\infty,a_1,\dots,a_{p^{d}}$. For any integers $i,j$ with $0\le i\leq j\le p^d$, we denote the subset $\{a_i,a_{i+1},\cdots,a_j\}$ by $D_{[i,j]}$ and for $1\leq w\leq p^d$, we set $\cD_w=\{D_{[i,w+i-1]}: 0\leq i\leq p^d-w+1\}$. Now we take $z_1+1$ distinct elements $A_0,\cdots,A_{z_1}$ from $\cD_w$, where $1\leq z_1\leq p^d-w+1$ and $A_i=D_{[k_i,w+k_i-1]}$. Without loss of generality, we assume that $k_i<k_j$ whenever $i<j$. Let
$$b_{i,j}=\left\{
\begin{aligned}
a_{k_i},\ \quad&\ \textnormal{for}\ i<j,   \\
a_{w+k_i-1},&\ \textnormal{for}\ i>j,
\end{aligned}
\right.$$
and $T_{i,j}$ be a subset of the system of coset representatives with respect to $L_{b_{i,j}}$ of size $z_2$ for $0\le i\ne j\le z_1$. Note that $1\leq z_2\leq p^{sd}$.
Then we define a family $\mathcal{S}=\{S_{i,j}\}_{0\le i,j\le z_1}$ of subsets of $G$ by
\begin{equation}\label{Eqn_psf0}
\begin{aligned}
&S_{i,i}=\bigcup_{a\in A_i}L_a\setminus\{0_G\}, \textup{ for }0\leq i\leq z_1\\
&S_{i,j}=\bigcup_{g\in T_{i,j}}(g+L_{b_{i,j}}), \textup{ for }0\le i\ne j\le z_1.
\end{aligned}
\end{equation}

We now claim that
{\lemma\label{lem 5}
 The family $\mathcal{S}$ defined by $(\ref{Eqn_psf0})$ is a partial sum family if and only if $z_2=w-1\ \textnormal{or}\ w$ for $2\le w\le p^d$, and $z_2=1$ for $w=1$. Furthermore, the parameters $(m,n,k,\lambda,\mu,t)=(z_{1}+1,p^{2sd},(p^{sd}-1)w+z_{1}z_{2}p^{sd},p^{sd}+w^2-3w+z_{1}z_{2}^2,w^2-w+z_{1}z_{2}^2,p^{sd}w-w+z_{1}z_{2}^2).$}
 \proof{By the definition of $S_{i,i}$'s, it is obvious that $0_G$ is not contained in each $S_{i,i}$. Since $|S_{i,i}|=(p^{sd}-1)w$ for any $i$ and $|S_{i,j}|=p^{sd}z_2$ for any $i\ne j$, we have $\sum_{j=0}^{z_1}|S_{i,j}|=\sum_{j=0}^{z_1}|S_{j,i}|=(p^{sd}-1)w+z_1z_2p^{sd}$, which is a constant integer for each $i$. Thus conditions (i) and (ii) for $\{S_{i,j}\}$ to form a partial sum family hold. What we need to show is that condition (iii) in Definition \ref{def_psf} holds for $\{S_{i,j}\}$ for three constants $\mu,\beta$ and $\gamma$ if and only if $z_2=w-1\ \textnormal{or}\ w$.

Rewrite condition (iii) in Definition \ref{def_psf} as the combination of two equations:
\begin{equation}\label{Eqn_psf1}
S_{i,i}^{2}+\sum_{l=0,l\ne i}^{z_1}S_{i,l}S_{l,i}=\mu G+\beta S_{i,i}+\gamma 0_G
\end{equation}
for any $0\leq i\leq z_1$, and
\begin{equation}\label{Eqn_psf2}
S_{i,j}(S_{i,i}+S_{j,j})+\sum_{l=0,l\ne i,j}^{z_1}S_{i,l}S_{l,i}=\mu G+\beta S_{i,j}
\end{equation}
for any $0\le i\ne j\le z_1$. By the definition of $S_{i,j}$'s, we have
$$S_{i,l}S_{l,j}=\sum_{g\in T_{i,l}}(g+L_{b_{i,l}})\sum_{h\in T_{l,j}}(h+L_{b_{l,j}})=\sum_{g\in T_{i,l}}\sum_{h\in T_{l,j}}(g+h+G)=|T_{i,l}||T_{l,j}|G=z_2^2G$$
for any $l\ne i,j$, as a direct result of $L_{b_{i,l}}L_{b_{l,j}}=G$. Thus we have
\begin{equation}\label{Eqn_psf3}
S_{i,l}S_{l,j}=z_2^2G.
\end{equation}
Plugging $(\ref{Eqn_psf3})$ into $(\ref{Eqn_psf1})$ and $(\ref{Eqn_psf2})$, we have
\begin{equation}\label{Eqn_psf4}
S_{i,i}^{2}=(\mu-z_1z_2^2)G+\beta S_{i,i}+\gamma 0_G
\end{equation}
for $0\leq i\leq z_1$, and
\begin{equation}\label{Eqn_psf5}
S_{i,j}(S_{i,i}+S_{j,j})=(\mu-(z_1-1)z_2^2)G+\beta S_{i,j}
\end{equation}
for $0\le i\ne j\le z_1$.
It suffices to show that (\ref{Eqn_psf4}) and (\ref{Eqn_psf5}) hold if and only if $z_2=w-1$ or $w$.

Assume that (\ref{Eqn_psf4}) and (\ref{Eqn_psf5}) hold. Then we will take characters on these two equations and get what we want.
By Lemma $\ref{lem_char}$, characters of $G$ are all given as
$$\chi_{b}:(x,y)\mapsto\chi_{b}(x,y)=\psi(b_{1}x+b_{2}y),$$
for $b=(b_{1},b_{2})\in G$. By taking the principal character (denoted simply by) $\chi_{0}$ on (\ref{Eqn_psf4}), we have
\begin{equation}\label{Eqn_psf6}
|S_{i,i}|^2=(\mu-z_1z_2^2)|G|+\beta|S_{i,i}|+\gamma,
\end{equation}
which implies that
\begin{equation}\label{Eqn_psf7}
(p^{sd}-1)^2w^2=(\mu-z_1z_2^2)p^{2sd}+\beta(p^{sd}-1)w+\gamma.
\end{equation}
\noindent By taking $\chi_{0}$ on $(\ref{Eqn_psf5})$, we have
\begin{equation}\label{Eqn_psf8}
2p^{sd}(p^{sd}-1)z_1z_2=(\mu-(z_1-1)z_2^2)p^{2sd}+\beta p^{sd}z_2.
\end{equation}
\noindent By taking a nonprincipal character $\chi_{b}$ on $L_{a}$ for $a\in J'$, we have
$$\chi_{b}(L_{a})=\left\{
\begin{aligned}
p^{sd},\ \textnormal{if}\ \chi_{b}&\ \textnormal{is principal on }L_{a},           \\
0,\quad\quad\quad&\textnormal{otherwise}.
\end{aligned}
\right.$$
As discussed in Section 2.1, if a character $\chi$ is principal on both $L_{a_{1}}$ and $L_{a_{2}}$ for $a_{1}\ne a_{2}\in J'$, then $\chi$ is principal on $G$. Thus we have
$$\chi_{b}(S_{i,i})=\left\{
\begin{aligned}
p^{sd}-w,\ &\textnormal{if}\ \chi_{b}\ \textnormal{is principal on some}\ L_{a}\ \textnormal{with}\ L_{a}\setminus\{0\}\subseteq\ S_{i,i},           \\
-w\quad,\ &\quad\quad\quad\quad\textnormal{otherwise}.
\end{aligned}
\right.$$
\\By taking a nontrivial character $\chi_{b}$ on $(\ref{Eqn_psf4})$, we have
\begin{equation}\label{Eqn_psf9}
\chi_{b}(S_{i,i})^{2}=\beta\chi_{b}(S_{i,i})+\gamma,
\end{equation}
As for $(\ref{Eqn_psf9})$, since for each nontrivial character $\chi_{b}$, $\chi_{b}(S_{i,i})\in\{p^{sd}-w,-w\}$, we have $p^{sd}-w+(-w)=\beta,$ and $(p^{sd}-w)(-w)=-\gamma$, then it implies that
$$\beta=p^{sd}-2w,\gamma=p^{sd}w-w^2.$$
\noindent From $(\ref{Eqn_psf7})$ we get the expression of $\mu$ in terms of $w,z_1$ and $z_2$,
$$\mu=w^2-w+z_1z_2^2,$$
and thus from $(\ref{Eqn_psf8})$ we deduce that
$$z_2^2-(2w-1)z_2+w^2-w=0,$$
i.e.,
$$z_2=w-1\ \textnormal{or}\ w.$$
Note that $1\leq w\leq p^d$ and $1\leq z_2\leq p^{sd}$, the integer $z_2$ can only take $w$ if $w=1$.

Conversely, suppose that $z_2=w$ for $w=1$, and $z_2=w-1$ or $w$ for $2\leq w\leq p^d$. By taking a nontrivial character $\chi_b$ on $(\ref{Eqn_psf5})$, we have
\begin{equation}\label{Eqn_psf10}
\chi_{b}(S_{i,j})\chi_{b}(S_{i,i}+S_{j,j})=\beta\chi_{b}(S_{i,j}).
\end{equation}
By Lemma $\ref{lem_4}$, it is sufficient to show that $(\ref{Eqn_psf7})$,\ $(\ref{Eqn_psf8})$,\ $(\ref{Eqn_psf9})$ and $(\ref{Eqn_psf10})$ hold for any nontrivial character $\chi_b$. Let $\mu=w^2-w+z_1z_2^2$, $\beta=p^{sd}-2w$ and $\gamma=p^{sd}w-w^2$, then $(\ref{Eqn_psf7})$ and $(\ref{Eqn_psf8})$ hold. $(\ref{Eqn_psf9})$ holds since $\chi_b(S_{i,i})\in\{p^{sd}-w,-w\}$ for any $0\le i\le z_1$ and nontrivial character $\chi_b$. Now we're left with the validity of $(\ref{Eqn_psf10})$. Since $\beta=p^{sd}-2w$ and
$$\chi_{b}(S_{i,i}+S_{j,j})\in\{-2w,p^{sd}-2w,2p^{sd}-2w\}$$
for any two distinct subsets $S_{i,i}$ and $S_{j,j}$ and nontrivial character $\chi_b$, we indicate that $\chi_{b}(S_{i,j})=0$ for any nontrivial character $\chi_{b}$ such that $\chi_{b}(S_{i,i}+S_{j,j})\ne\beta$. Let $\chi_{b}$ be such a character, then $\chi_{b}(S_{i,i})=\chi_{b}(S_{j,j})=p^{sd}-w$ or $-w$. In both cases we can deduce that $\chi_{b}$ is nonprincipal on $L_{b_{i,j}}$, and thus
$$\chi_{b}(S_{i,j})=\chi_{b}(\bigcup_{g\in T_{ij}}(g+L_{b_{i,j}}))=\chi_{b}(L_{b_{i,j}})\sum_{g\in T_{i,j}}\chi_{b}(g)=0.$$
Hence $(\ref{Eqn_psf10})$ holds.\qed
}
We can now prove Theorem \ref{thm_psfmain}.

\proof[Proof of Theorem \ref{thm_psfmain}]{
By Lemma \ref{lem_DF} and Lemma \ref{lem 5}, we attain a DSRG with parameters mentioned in this theorem.
}

\begin{remark}
Denote the digraph associated with $\mathcal{S}$ by $\mathcal{G}$. Since $w\le p^d$, the sets $S_{i,i}$'s for $0\leq i\leq z_1$ are all not equal to $G\setminus\{0_G\}$, which means that $\mathcal{G}$ is not complete. We need to rule out the cases $t=k$ and $t=0$ which also make $\mathcal{G}$ trivial. It is easy to see that $t=p^{sd}w-w+z_1z_2^2>0$.  When the case $t=k$ happens, we get $z_2=p^{sd}$. Eventually, we conclude that $\mathcal{G}$ is a (non-complete) SRG if and only if $w=z_2=p^d$ and $s=1$. By the way, since in the case that $s=1$, the local ring is isomorphic to a finite field, we call this case is the ``finite field" case.
\end{remark}

We denote the expressions of $k,\lambda,\mu,t$ in terms of $w,z_1$ and $z_2$ in the above theorem by $k(w,z_1,z_2),\lambda(w,z_1,z_2),\mu(w,z_1,z_2),t(w,z_1,z_2)$ for convenience, respectively.

\begin{remark}
The 16 new tuples of parameters listed in Section 1 are pairwise complementary and can all be obtained from the ``finite field" case. If a $(z_1+1,p^{2d},k(w,z_1,z_2),\lambda(w,z_1,\\z_2),\mu(w,z_1,z_2),t(w,z_1,z_2))$-PSF comes from our construction with $s=1$, then so does another one with complementary parameters $(z_1+1,p^{2d},k(p^d+1-w,z_1,p^d-z_2),\lambda(p^d+1-w,z_1,p^d-z_2),\mu(p^d+1-w,z_1,p^d-z_2),t(p^d+1-w,z_1,p^d-z_2))$.
\end{remark}

Let $\mathcal{S}=\{S_{i,j}\}_{0\leq i,j\leq z_1}$ be a family obtained from $(\ref{Eqn_psf0})$, and $c$ be a fixed integer with $0\leq c\leq z_1$. Then the deleted family $\mathcal{S}'=\{S_{i,j}\}_{0\leq i,j\leq z_1,i,j\ne c}$ of $\mathcal{S}$ is a PSF with parameters $(z_1,p^{2sd},k(w,z_1-1,z_2),\lambda(w,z_1-1,z_2),\mu(w,z_1-1,z_2),t(w,z_1-1,z_2))$. Furthermore, in the ``finite field" case, we can choose the subsets $A_0,A_1,\cdots,A_{z_1}$ appropriately so that the resulting partial sum family is uniform.

\begin{corollary}
For any prime integer $p$ and positive integers $d,w$ such that $w|p^d+1$, there exists a uniform partial sum family $\{S_{i,j}\}_{0\leq i,j\leq z_1}$ with parameters $(z_1,p^{2d},k(w,z_1-1,z_2),\lambda(w,z_1-1,z_2),\mu(w,z_1-1,z_2),t(w,z_1-1,z_2))$, where $z_1=\frac{p^d+1}{w}$ and $z_2=w-1$ or $w$.
\end{corollary}
\proof{
Set $z_1=\frac{p^d+1}{w}$ and $A_i=T_{[wi,w(i+1)-1]}$ for $0\leq i\leq z_1$. Let $\{S_{i,j}\}_{0\leq i,j\leq z_1}$ be the family defined in (\ref{Eqn_psf0}) by $\{A_0,\cdots, A_{z_1}\}$ with $z_2=w-1$ or $w$ . As shown in Lemma \ref{lem 5}, this family is a partial sum family with parameters $(z_1,p^{2d},k(w,z_1-1,z_2),\lambda(w,z_1-1,z_2),\mu(w,z_1-1,z_2),t(w,z_1-1,z_2))$. It is easy to see that $\{S_{i,i}:0\leq i\leq z_1\}$ forms a partition of $G\setminus \{0\}$. As we discussed in the proof of Lemma \ref{lem 5}, $|S_{i,i}|=(p^{d}-1)w$ and $|S_{i,j}|=p^{d}z_2$ for any $i\neq j$ which ensures conditions (i) and (ii) in Definition \ref{def_upsf} to be satisfied. Hence this partial sum family is uniform.\qed
}

\section{Constructions of Directed Strongly Regular Graphs from Semi-direct Products of Groups}
In this section we will use the notations mentioned in Subsection 2.3. Denote the additive group of $\mathbb{F}_{q}$ by $H$. Now let $\rho$ be an automorphism of $H$ such that $\rho(x)=\omega^e x$ for $x\in\mathbb{F}_{q}$, then $K=\langle\rho\rangle$ is an automorphism group of $H$ which is isomorphic to $C_{0}$. Let $G=H\rtimes K$. For any subset $S$ of $\mathbb{F}_{q}$, we say $S$ is $K$-invariant if $\phi(S)=S$ holds for any $\phi\in K$. Note that $S$ is $K$-invariant if and only if $S$ is the union of some cyclotomic cosets of $C_{0}$.
\\Suppose $D$ is a $K$-invariant subset of $H$. Denote $D_1=-1+D=\{-1+d|d\in D\}$ and let $W=(D_1,K)\subset G$. We compute in $\mathbb{Z}[H]$ that
$$\begin{aligned}
W^2&=(D_1,K)(D_1,K)\\
&=\sum_{d\in D}\sum_{i=0}^{f-1}(-1+d,\rho^i)\sum_{d'\in D}\sum_{j=0}^{f-1}(-1+d',\rho^j)\\
&=\sum_{d\in D}\sum_{d'\in D}\sum_{i=0}^{f-1}\sum_{j=0}^{f-1}(-1+d+\omega^{ei}(-1+d'),\rho^{i+j}).
\end{aligned}$$
Substituting $i+j$ with $j$, we have
$$\begin{aligned}
W^2&=\sum_{d\in D}\sum_{d'\in D}\sum_{i=0}^{f-1}\sum_{j=0}^{f-1}(-1+d-\omega^{ei}+\omega^{ei}d'),\rho^{j})\\
&=\sum_{j=0}^{f-1}\sum_{d\in D}\sum_{i=0}^{f-1}\sum_{d'\in D}(-1+d-\omega^{ei}+\omega^{ei}d'),\rho^{j}).
\end{aligned}$$
And since $D$ is $K$-invariant, we deduce that
$$\begin{aligned}
W^2&=\sum_{j=0}^{f-1}\sum_{d\in D}\sum_{i=0}^{f-1}\sum_{d'\in D}(-1+d-\omega^{ei}+d'),\rho^{j})\\
&=((-1+D)C_0^{(-1)}D,K)\\
&=(D_1DC_0^{(-1)},K).
\end{aligned}$$

Then we have the following result:
\begin{lemma}
If $DC_0^{(-1)}=a[0]+bH$ for integers $a$ and $b$, then $W$ is the connection set of a directed strongly regular Cayley graph.
\end{lemma}
\proof{
Since $DC_0^{(-1)}=a[0]+bH$, we have $W^2=(D_1DC_0^{(-1)},H)=(D_1(a[0]+bH),K)=(aD_1+b|D_1|H,K)=aW+b|D|G$, then by Lemma 2.5 it completes the proof.\qed
}
In the case that $C_0$ is not contained in $D$, we deduce that $a=-b$ by comparing the coefficients of [0] in both sides of $DC_0^{(-1)}=a[0]+bH$. Then $b=\frac{|C_0||D|}{q-1}$ and $DC_0^{(-1)}=\frac{|C_0||D|}{q-1}H^{\ast}$, where $H^{\ast}=H\setminus\{0\}$.
\\In the rest of this section, our task is to find the subset $D$ such that $DC_0^{(-1)}=a[0]+bH$.
If we require $f$ is even, then since $-1\in C_0$ we have $C_0^{(-1)}=C_0$ and by setting $D$ to be some specific cyclotomic coset we can get a DSRG. Suppose $f$ is even and $D=C_{i}$ for some $1\le i\le f-1$, then by (\ref{Eqn_schurring}), we have
$$DC_0^{(-1)}=C_iC_0=\sum_{k=0}^{e-1}(i,k)C_{k}+|C_{0}\cap-C_{i}|[0]=\sum_{k=0}^{e-1}(i,k)C_{k}.$$
Then we have
\begin{lemma}
Suppose that for some $1\le i\le e-1$, $(i,j)$'s are all the same for any $0\le j\le f-1$. Let $D=C_i$, then $W$ is the connection set of a directed strongly regular Cayley graph.
\end{lemma}
\proof{
By assumption, for any $j$, $(i,j)$ is a constant integer, which we say $s$. Then $$DC_0^{(-1)}=\sum_{k=0}^{e-1}(i,k)C_{k}=\sum_{k=0}^{e-1}sC_{k}=sH^{\ast},$$
and then by Lemma 4.1 it completes the proof.\qed
}
\begin{lemma}
(\cite{9})
When $e=4$ and $f$ is even, the cyclotomic matrix
$$[(i,j)]=
\begin{bmatrix}
(0,0)  & (0,1)  & (0,2)  & (0,3)\\
(0,1)  & (0,3)  & (1,2)  & (1,2)\\
(0,2)  & (1,2)  & (0,2)  & (1,2)\\
(0,3)  & (1,2)  & (1,2)  & (0,1)\\
\end{bmatrix},$$
where
$$(0,0)=\frac{q-11-6s}{16},\ (0,1)=\frac{q-3+2s+8t}{16},\ (0,2)=\frac{q-3+2s}{16},$$ $$(0,3)=\frac{q-3+2s-8t}{16},\ (1,2)=\frac{q+1-2s}{16},$$
and $s,t$ is such that $q=s^2+4t^2$, $s\equiv1$ (mod $4$) and gcd$(p,s)=1$ provided that $p\equiv1$ (mod $4$).
\end{lemma}
By Lemma 4.2, we choose $D=C_2$ and require $(0,2)$ and $(1,2)$ appearing in Lemma 4.3 to be equal, which follows that
$$\frac{q-3+2s}{16}=\frac{q+1-2s}{16},$$
that is, $s=1$. Therefore $q$ is a prime power of form $1+4t^2$. By \cite{9}, we know that the Diophantine equation
$$x^n-y^2=1$$
has no integer solutions for any $n>1$. Then $q$ must be an odd prime, i.e., $q=p$, and we obtain the following result.
\begin{thm}
Let $p=1+4t^2$ be a prime with $t$ even, then $W=(-1+C_2,H)$ is the connection set of a directed strongly regular graph with parameters $$(\frac{(p-1)p}{4},\frac{(p-1)^2}{16},\frac{(p-1)^2}{64},\frac{(p-1)^2}{64}-\frac{p-1}{16},\frac{(p-1)^2}{64}).$$
\end{thm}
There are 59 such primes $p$ less that $10^6$: 17, 257, 401, 577, 1297, 1601, 3137, 7057, 13457, 14401,
15377, 24337, 25601, 30977, 32401, 33857, 41617, 50177, 55697, 57601, 65537, 67601, 69697, 78401,
80657, 90001, 115601, 147457, 156817, 160001, 176401, 190097, 193601, 197137, 215297, 246017,
287297, 295937, 309137, 331777, 341057, 404497, 414737, 462401, 484417, 490001, 495617, 512657,
547601, 577601, 583697, 608401, 614657, 665857, 739601, 746497, 846401, 876097, 921601.

Similarly, we treat the case that $e=6$. In \cite{10}, the cyclotomic numbers of order 6 for a finite field $\mathbb{F}_q$ are all available. In this paper we need the ones for the case that 2 is a cube in $\mathbb{F}_q$.
\begin{lemma}
(\cite{10})
Let $e=6$, $f$ even and $m\equiv0$ (mod $3$), where the integer $m$ is such that $\omega^m=2$. Then the cyclotomic matrix
$$[(i,j)]=
\begin{bmatrix}
(0,0)  & (0,1)  & (0,2)  & (0,3)  & (0,4)  & (0,5)\\
(0,1)  & (0,5)  & (1,2)  & (1,3)  & (1,4)  & (1,2)\\
(0,2)  & (1,2)  & (0,4)  & (1,4)  & (2,4)  & (1,3)\\
(0,3)  & (1,3)  & (1,4)  & (0,3)  & (1,3)  & (1,4)\\
(0,4)  & (1,4)  & (2,4)  & (1,3)  & (0,2)  & (1,2)\\
(0,5)  & (1,2)  & (1,3)  & (1,4)  & (1,2)  & (0,1)\\
\end{bmatrix}.$$
where the possible distinct cyclotomic numbers in the cyclotomic matrix are expressed as follows,
$$(0,0)=\frac{q-17-20a}{36},\ (0,1)=\frac{q-5+4a+18b}{36},\ (0,2)=\frac{q-5+4a+6b}{36},\ (0,3)=\frac{q-5+4a}{36},$$ $$(0,4)=\frac{q-5+4a-6b}{36},\ (0,5)=\frac{q-5+4a-18b}{36},\ (1,2)=(1,3)=(1,4)=(2,4)=\frac{q+1-2a}{36},$$
and $a,b$ is such that
$$q=a^2+3b^2,\ a\equiv1\ (\textnormal{mod }3),\ \textnormal{gcd}(p,a)=1,$$
and
$$\omega^{(q-1)/3}\equiv-(a+b)/(a-b)\ (\textnormal{mod }p).$$
\end{lemma}
Applying Lemma 4.2 again, we choose $D=C_3$ and require $(0,3)_6=(1,3)_6$, which leads to
$$\frac{q-5+4a}{16}=\frac{q+1-2a}{16},$$
i.e., $a=1$. Therefore $q$ is a prime power of form $1+3b^2$.
\begin{thm}
Suppose that $q=1+3b^2$ is a prime power with $b$ even, and that $2$ is a cube in $\mathbb{F}_{q}$, then $W=(-1+C_3,H)$ is the connection set of a directed strongly regular graph with parameters $$(\frac{(q-1)q}{6},\frac{(q-1)^2}{36},\frac{(q-1)^2}{216},\frac{(q-1)^2}{216}-\frac{q-1}{36},\frac{(q-1)^2}{216}).$$
\end{thm}
There are 25 of such $q$ less than $10^6$, all of which are primes: 109, 433, 3889, 18253, 21169, 43201, 47629, 73009,84673, 139969, 147853, 172801, 181549, 209089, 238573, 259309, 270001, 280909, 326701, 363313, 415153, 544429, 640333, 762049, 995329.

\section{Concluding Remark}
One of the main results in this paper is to construct partial sum families with local rings. From this construction we obtain a new infinite family of DSRGs, which are directed strongly regular $m$-Cayley digraphs for each $m\geq 2$. Some of these DSRGs have new parameters, which implies that PSFs are useful in obtaining new DSRGs. Corollary 3.5 tells that uniform PSFs defined in \cite{7} are available from our first construction. We also obtain two infinite families of directed strongly regular Cayley graphs from semi-direct products of groups. These two constructions do not produce DSRGs with new feasible parameters, but they produce DSRGs from new methods in non-abelian groups.

\section*{Acknowledgement}
The work of the first two authors was supported by National Natural Science Foundation of China under Grant No. 11771392.

\end{document}